\documentclass[11pt] {article}
\usepackage{authblk}
\usepackage{etex}
\usepackage[utf8]{inputenc}
\usepackage[english]{babel}
\usepackage{amsmath}
\usepackage{amsthm}
\usepackage{amssymb}
\usepackage{sectsty}
\usepackage{titlesec}
\usepackage{color}
\usepackage[color,matrix,arrow]{xy}
\usepackage{amsgen}
\usepackage{amstext}
\usepackage{amsbsy}
\usepackage{amsopn}
\usepackage{amsfonts}
\usepackage{eepic}
\usepackage{graphicx}
\usepackage{epsf}
\usepackage{mathtools}
\usepackage{pstricks}
\usepackage{float}
\usepackage{enumerate}
\usepackage{tikz}
\usepackage{faktor}
\usepackage{tikz-cd}
\usetikzlibrary{decorations.pathmorphing}
\xyoption{all}
\usetikzlibrary{patterns}
 
\usepackage{pgf}
\usepackage{tikz}

\newcommand{\footrecall}[1]{}

\titleformat*{\section}{\large\bfseries}
\titleformat*{\subsection}{\normalsize \bfseries}

\title{On uniformly continuous endomorphisms of hyperbolic groups}

\author{Andr\'e Carvalho \thanks{andrecruzcarvalho@gmail.com}}
\affil{Centre of Mathematics, University of Porto, R. Campo Alegre, 4169-007 Porto,
Portugal}

\date{}

\begin{document}

\newtheorem{theorem}{Theorem}[section]
\newtheorem{lemma}[theorem]{Lemma}
\newtheorem{question}[theorem]{Question}
\newtheorem{remark}[theorem]{Remark}
\newtheorem{proposition}[theorem]{Proposition}
\newtheorem{corollary}[theorem]{Corollary}

\theoremstyle{definition}
\newtheorem{definition}[theorem]{Definition}
\newtheorem{example}[theorem]{Example}

\newcommand{\ophi}{\overline{\varphi}}
\newcommand{\opsi}{\overline{\psi}}
\newcommand{\N}{\mathbb{N}}
\newcommand{\Z}{\mathbb{Z}}
\newcommand{\F}{\mathbb{F}}
\newcommand{\R}{\mathbb{R}}
\newcommand{\Q}{\mathbb{Q}}
\newcommand{\Fix}{\text{Fix}}
\newcommand{\Sing}{\text{Sing}}
\newcommand{\Aut}{\text{Aut}}
\newcommand{\Per}{\text{Per}}
\newcommand{\Reg}{\text{Reg}}
\newcommand{\End}{\text{End}}
\newcommand{\mc}{\mathcal}

\maketitle

\begin{abstract}
We prove a generalization of the fellow traveller property for a certain type of quasi-geodesics and use it to present three equivalent geometric formulations of the bounded reduction property and prove that it is equivalent to preservation of a coarse median. We then provide an affirmative answer to a question from Ara\'ujo and Silva as to whether every nontrivial uniformly continuous endomorphism of a hyperbolic group with respect to a visual metric satisfies a H\"older condition. We remark that these results combined with the work done by Paulin prove that every endomorphism admitting a continuous extension to the completion has a finitely generated fixed point subgroup. 
\end{abstract}

\section{Introduction}
\label{intro}
The dynamical study of endomorphisms of groups started with the (independent) work of Gersten \cite{[Ger87]} and Cooper \cite{[Coo87]}, using respectively graph-theoretic and topological approaches. They proved that the subgroup of fixed points $\Fix(\varphi)$ of some fixed automorphism $\varphi$ of $F_n$ is always finitely generated, and Cooper succeeded on classifying from the dynamical viewpoint the fixed points of the continuous extension of $\varphi$ to the boundary of $F_n$. Bestvina and Handel subsequently developed the theory of train tracks to prove that $\Fix(\varphi)$ has rank at most $n$ in \cite{[BH92]}. The problem of computing a basis for $\Fix(\varphi)$ had a tribulated history and was finally settled by Bogopolski and Maslakova in 2016 in \cite{[BM16]}.

This line of research extended early to wider classes of groups. For instance, Paulin proved in 1989 that the subgroup of fixed points of an automorphism of a hyperbolic group is finitely generated \cite{[Pau89]}. Fixed points were also studied for right-angled Artin groups \cite{[RSS13]} and lamplighter groups \cite{[MS18]}.

Regarding the continuous extension of an endomorphism to the completion when a suitable metric is considered 
, infinite fixed points of automorphisms 
of free groups were also discussed by Bestvina and Handel in  \cite{[BH92]} and Gaboriau, Jaeger, Levitt and Lustig in \cite{[GJLL98]}. The dynamics of free groups automorphisms is proved to be asymptotically periodic in \cite{[LL08]}. In \cite{[CS09a]}, the dynamical study of infinite fixed points was performed for monoids defined by special confluent rewriting systems (which contain free groups as a particular case). This was also achieved in \cite{[Sil13]} for virtually injective endomorphisms of virtually free groups.
Endomorphisms of free-abelian times free groups $\Z^m\times F_n$ have been studied in \cite{[DV13]} and their continuous extension to the completion 
 is studied in \cite{[Car20]} and the case of $F_n\times F_m$ with $m,n\geq 2$ is dealt with in \cite{[Car20b]}
 
For natural reasons, in order to study the dynamics of the continuous extension to the completion, in case the topology is defined by a given distance (as it is in the case of hyperbolic groups via visual metrics), it is of utmost importance to describe uniformly continuous endomorphisms for a certain class of groups. Also, it is important to notice that one of the essential tools used in proving these results is the bounded reduction property (also known as the bounded cancellation lemma) introduced in \cite{[Coo87]} and followed by many others.

Coarse median preservation was recently introduced by Fioravanti and it turns out to be a useful tool to obtain interesting properties of automorphisms (see \cite{[Fio21]}), including finiteness results on fixed subgroups. We will prove that coarse-median preserving endomorphisms are precisely the ones for which the BRP holds.

Motivated by the possibility of defining new pseudometrics in the group of automorphisms of a hyperbolic group, the authors in \cite{[AS16]} studied endomorphisms for which a H\"older condition holds.

In this paper, we will describe uniformly continuous endomorphisms of hyperbolic groups endowed with appropriate distances and present some geometric versions of the bounded reduction property with hope that these techniques could be used to study the dynamics of points at infinity
for arbitrary uniformly continuous endomorphisms of hyperbolic groups, for which not much is known so far. We will also show that H\"older conditions are satisfied by every uniformly continuous endomorphism with respect to a visual metric, answering a question from \cite{[AS16]} and highlight that this result, combined with previous work from Paulin suffices to show that every uniformly continuous endomorphism of a hyperbolic group has finitely generated fixed point subgroup.

The paper is organized as follows. In Section \ref{secpreliminares}, we present some preliminaries on hyperbolic metric spaces and hyperbolic groups. In Section \ref{secftp}, we present a generalization of the fellow traveller property for $(1,k)$-quasi-geodesics, which will be seen to arise naturally. We formulate the bounded reduction property in geometric terms throughout Section \ref{secbrp} and prove that it coincides with coarse-median preservation. In Section \ref{secuc}, we show that every nontrivial uniformly continuous endomorphism of a hyperbolic group must satisfy a H\"older condition, answering a question by Ara\'ujo and Silva, and observe that these results combined with the work in \cite{[Pau89]} yield that every endomorphism of a finitely generated hyperbolic group admitting a continuous extension to the completion has a finitely generated fixed point subgroup.

\section{Preliminaries}
 \label{secpreliminares}
We now introduce some well-known results on hyperbolic groups. For more details, the reader is referred to \cite{[GH90]} and \cite{[BH99]}.

A mapping $\varphi:(X,d)\to(X',d')$ between metric spaces is called an \emph{isometric embedding} if $d'(x\varphi,y\varphi)=d(x,y)$, for all $x,y\in X$. A surjective isometric embedding is an \emph{isometry}.

A metric space $(X, d)$ is said to be \emph{geodesic} if, for all $ x, y \in X$, there exists an isometric embedding
$\xi: [0, s]\to X$ such that $ 0\xi = x$ and $s\xi = y$, where $[0, s] \subset \R$ is endowed with the usual metric of $\R$.
We call $\xi$ a geodesic of $(X, d)$. We shall often call Im$(\xi)$ a geodesic as well. In this second sense, we
may use the notation $[x, y]$ to denote an arbitrary geodesic connecting $x$ to $y$.  When the endpoint of a geodesic $\alpha$ coincides with the starting point of a geodesic $\beta$, we denote the concatenation of both geodesics by $\alpha+\beta$. Note that a geodesic
metric space is always (path) connected. A {\em quasi-isometric embedding} of metric spaces is a mapping $\varphi:(X,d) \to
(X',d')$ such that there exist constants $\lambda \geq 1$ and $K \geq 0$
satisfying
$$\frac{1}{\lambda}d(x,y) -K \leq d'(x\varphi,y\varphi) \leq \lambda d(x,y) +K$$
for all $x,y \in X$. We may call it a
$(\lambda,K)$-quasi-isometric embedding if we want to stress the
constants. 

If in addition
$$\forall x'\in X\,\exists x\in X: d'(x',x\varphi)\leq K,$$
we say that $\varphi$ is a \emph{quasi-isometry}.
Two metric spaces $(X, d)$ and $(X', d')$ are said to be \emph{quasi-isometric} if there exists a quasi-isometry
$\varphi : (X, d) \to (X', d')$. Quasi-isometry turns out to be an equivalence relation on the class of metric
spaces.
 A \emph{$(\lambda,K)$-quasi-geodesic} of $(X, d)$ is a $(\lambda,K)$-quasi-isometric embedding $\xi : [0, s] \to X$ such that $0\xi = x$
and $s\xi = y$, where $[0, s] \subset \R$ is endowed with the usual metric of $\R$.

Given $x,y,z\in X$, a \emph{geodesic triangle} $[[x,y,z]]$ is a collection of three geodesics $[x,y]$, $[y,z]$ and $[z,x]$ in $X$. Given $\delta\geq 0$, we say that $X$ is \emph{$\delta$-hyperbolic} if $$\forall w\in [x,y]\, d(w,[y,z]\cup [z,x])\leq \delta$$ holds for every geodesic triangle $[[x,y,z]]$.

Let $(X,d)$ be a metric space and $Y,Z$ nonempty subsets of $X$. We call the $\varepsilon$-neighbourhood of $Y$ in $X$ and we denote by $\mc V_\varepsilon(Y)$ the set $\{x\in X\mid d(x,Y)\leq \varepsilon\}$.
We call the {\em Hausdorff distance} between
$Y$ and $Z$ and we denote by Haus$(Y,Z)$, the number defined by
$$\inf\{ \varepsilon>0\mid Y\subseteq \mc V_\varepsilon (Z) \text{ and } Z\subseteq\mc V_\varepsilon(Y)\}$$
if it exists. If it doesn't, we say that Haus$(Y,Z)=\infty$.

Given a group $H=\langle A\rangle$, consider its Cayley graph $\Gamma_A(H)$ with respect to $A$ endowed with the \emph{geodesic metric} $d_A$, defined by letting $d_A(x,y)$ to be the length of the shortest path in $\Gamma_A(H)$ connecting $x$ to $y$. This is not a geodesic metric space, since $d_A$ only takes integral values. However,  we can define the \emph{geometric realization} $\bar\Gamma_A(H)$ of its Cayley graph $\Gamma_A(H)$ by embedding $(H,d_A)$ isometrically into it. Then, edges of the Cayley graph become segments of length $1$. With the metric induced by $d_A$, which we will also denote by $d_A$, $\bar\Gamma_A(H)$ becomes a geodesic metric space.

We say that a group $H$ is \emph{hyperbolic} if the metric space  $(\bar\Gamma_A(H),d_A)$ is hyperbolic. We will simply write $d$ instead of $d_A$ when no confusion arises. Also, for $x\in H$ we will often denote $d_A(1,u)$ by $|u|$.

Given a finite alphabet $A$, we write $\widetilde{A}=A\cup A^{-1}$ and $\widetilde{A}^*$ is the free monoid on $\widetilde{A}$. From now on, $H$ will denote a finitely generated hyperbolic group generated by a finite set $A$ and $\pi:\widetilde{A}^*\to H$ will be a matched epimorphism. A homomorphism $\pi:\widetilde{A}^*\to H$ is said to be \emph{matched} if $a^{-1}\pi=(a\pi)^{-1}$.

An important property of the class of automatic groups, for which the class of hyperbolic groups is a subclass, is the \emph{fellow traveller property}.  Given a word $u\in \widetilde{A}^*$, we denote by $u^{[n]}$, the prefix of $u$ with $n$ letters. If $n>|u|$, then we consider $u^{[n]}=u$. We say that the fellow traveller property holds for $L\subseteq \widetilde A^*$ if, for every $M\in \N$, there is some $N\in \N$ such that, for every $u,v\in L,$ $$d_A(u\pi,v\pi)\leq M \Rightarrow d_A(u^{[n]}\pi,v^{[n]}\pi)\leq N,$$
for every $n\in \N$.

Given $g,h,p\in H$, we define the \emph{Gromov product of $g$ and $h$} taking $p$ as basepoint by
$$(g|h)_p^A=\frac 12 (d_A(p,g)+d_A(p,h)-d_A(g,h)).$$ We will often write $(g|h)$ to denote $(g|h)_1^A$. Notice that, in the free group case, we have that $(g|h)=|g\wedge h|,$ where $u\wedge v$ denotes the longest common prefix between $u$ and $v$.

Let $G$ be a subgroup of a hyperbolic group $H = \langle A\rangle$ and
let $q \geq 0$. We
say that $G$ is $q$-{\em quasiconvex} with respect to $A$ if
$$\forall x \in [g,g']\hspace{.3cm} d_A(x,G) \leq q$$
holds for every geodesic $[g,g']$ in $\overline{\Gamma}_A(H)$ with endpoints in
$G$. We say that $G$ is {\em quasiconvex} if it is $q$-quasiconvex
for some $q \geq 0$.

There are many ways to describe the \emph{Gromov boundary} $\partial H$ of $H$, such as being the equivalence classes of geodesic rays, when two rays are considered equivalent if the Hausdorff distance between them is finite. 
Another model for $\partial H$ can be defined using \emph{Gromov sequences}. We say that a sequence of points $(x_i)_{i\in \N}$ in $H$ is a Gromov sequence if $(x_i|x_j)\to \infty$ as $i\to\infty$ and $j\to\infty.$ Two such sequences $(x_i)_{i\in \N}$ and $(y_j)_{j\in \N}$ are \emph{equivalent} if $$\lim\limits_{i,j\to \infty} (x_i|y_j)=\infty.$$ 
The set of all the equivalence classes is another standard model for the boundary $\partial H$. Identifying an element $h$ in $H$ with the constant sequence $(h)_i$, we can extend the Gromov product to the boundary by putting, for all $\alpha,\beta\in \partial H$,

$$(\alpha|\beta)_p^A=\sup\left\{\liminf_{i,j\to\infty} (x_i|y_j)_p^A \,\big\lvert\, (x_i)_{i\in \N}\in \alpha, (y_j)_{j\in \N}\in \beta\right\},
$$

We define 
$$\rho^A_{p,\gamma}(g,h) =
\left\{
\begin{array}{ll}
e^{-\gamma(g|h)_p^A}&\mbox{ if } g\neq h\\
0&\mbox{ otherwise}
\end{array}
\right.$$
for all $p,g,h \in H$.

Given $p \in H$, $\gamma > 0$ and $T \geq 1$, we
denote by $V^A(p,\gamma,T)$ the set of all metrics $d$ on $H$ such
that
\begin{align}
\label{ineq}
\frac{1}{T}\rho^A_{p,\gamma}(g,h) \leq d(g,h) \leq
T\rho^A_{p,\gamma}(g,h)
\end{align}

We refer to the metrics in some $V^A(p,\gamma,T)$ as the \emph{visual metrics} on $H$. Let $d\in V^A(p,\gamma,T)$. The metric space $(H,d)$ is not complete in general. However, its completion $(\widehat H,\hat d)$ is a compact space and can be obtained by considering $\widehat H= H\cup \partial H$. It is a well-known fact that the topology induced by $\hat d$ on $\partial H$ is the \emph{Gromov topology} and that all visual metrics originate equivalent completions. Therefore, we will write \emph{the completion} to mean a completion of the space when a visual metric is considered.

Considering the extension of $\rho^A_{p,\gamma}$ to the boundary, we define, for $\alpha,\beta\in \hat H$
$$\hat\rho^A_{p,\gamma}(\alpha,\beta) =
\left\{
\begin{array}{ll}
e^{-\gamma(\alpha|\beta)_p^A}&\mbox{ if } \alpha\neq \beta\\
0&\mbox{ otherwise}
\end{array}
\right.$$

By continuity, for every $\alpha,\beta\in \hat H$, the inequalities 
\begin{align*}
\frac{1}{T}\hat\rho^A_{p,\gamma}(\alpha,\beta) \leq \hat d(\alpha,\beta) \leq
T\hat\rho^A_{p,\gamma}(\alpha,\beta)
\end{align*}
hold  \cite[Section III.H.3]{[BH99]}.

\section{Fellow Traveller Property for $(1,k)$-quasi-geodesics}
 \label{secftp}
The goal of this section is to present a slightly more general formulation of the fellow traveller property. The result is most likely known but since it is not easy to find a proof of the result, we include one here, for sake of completeness. While the property is usually stated considering geodesics, we can see that we get essentially the same result when the paths considered are $(1,k)$-quasi-geodesics for some $k\in \N$. We will present the result in a different version, considering quasi-geodesics with endpoints at distance at most one from one another, but we remark that the result holds as long as the distance between the endpoints is bounded by some constant.

\begin{proposition}
\label{ftrav}
Let $u$ be a $(1,r)$-quasi-geodesic and $v$ be a $(1,s)$-quasi-geodesic with the same starting point. Then, there is a constant $N$ depending on $r,s,\delta$ such that,  for all $n\in \N$, we have that $$d_A(u\pi,v\pi)\leq 1\Rightarrow d_A(u^{[n]}\pi,v^{[n]}\pi)\leq N.$$
\end{proposition}
\noindent\textit{Proof.} Since the Cayley graph of $H$ with respect to $A$ is vertex-transitive, we can assume that $1$ is the starting point of both $u$ and $v$. Let $p$ and $q$ be the endpoints of $u$ and $v$, and consider a geodesic $w$ from $q$ to $p$. 
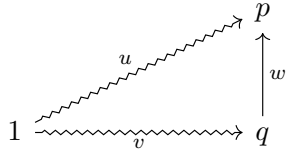
$$
\begin{tikzcd}[sep=large]
                                                                                                                               &            & p \\
   1 \ar[urr,squiggly,"u"] \ar[rr,swap,squiggly,"v"] &            & q  \ar[u,swap,"w"] 
\end{tikzcd}
$$
Let $k=\max\{r,s\}$. Then $u,v$ and $w$ are all $(1,k)$-quasi-geodesics. By Corollary 1.8, Chapter III.H in \cite{[BH99]}, there is a constant $\delta'$, depending only on $r,s,\delta$, such that the triangle $(u,v,w)$ is $\delta'$-thin. 
Let $n\in \N$.
We will prove that $$d_A(u^{[n]}\pi,v^{[n]}\pi)\leq 3k+2\delta'+4.$$
Consider the factorizations of $u$ and $v$ given by
$$\begin{tikzcd}
1 \ar[rr, squiggly,"u^{[n]}"] && p_n \ar[rr,squiggly,"u_n"] && p, 
\end{tikzcd}\quad 
\begin{tikzcd}
1 \ar[rr, squiggly,"v^{[n]}"] && q_n \ar[rr,squiggly,"v_n"] && q, 
\end{tikzcd}.$$
Suppose that $v_n=1$. This means that $q_n=q$, $d_A(1,q)\leq n$ and $d_A(1,p)\leq n+1.$ Also, notice that, since $u$ is a $(1,k)$-quasi-geodesic, we have that 
\begin{align*}
|u|-k\leq d_A(1,p)\leq |u|. 
\end{align*} 
We have that $$d_A(u^{[n]}\pi,p)\leq |u|-n\leq d_A(1,p)-n+k\leq k+1\leq 3k+2\delta'+4.$$
The case $u_n=1$ is analogous, so we assume that $u_n,v_n\neq 1$.

Since $(u,v,w)$ is a $\delta'$-thin triangle, there is some $m\leq |u|,$ such that $$d_A(v^{[n]}\pi,u^{[m]}\pi)\leq \delta'+1.$$ 
Again, since $u$ and $v$ are $(1,k)$-quasi-geodesics, we have that 
\begin{align*}
d_A(1,q)-n-k\leq |v|-n-k\leq d_A(q_n,q) \leq |v|-n\leq d_A(1,q)-n+k
\end{align*}
and
\begin{align*}
 d_A(1,p)-m-k\leq |u|-m-k\leq d_A(p_m,p)\leq |u|-m\leq d_A(1,p)-m+k.
 \end{align*}

So, 
\begin{align*}
d_A(1,q)-n-k&\leq d(q_n,q)\leq d(q_n,p_m)+d(p_m,p)+d(p,q)\\
&\leq \delta'+1+d_A(1,p)-m+k +1,
\end{align*}
thus, $m-n\leq \delta'+1+d_A(1,p)-d_A(1,q)+2k+1\leq \delta'+2k+3.$
Similarly, we have that 
\begin{align*}
d_A(1,p)-m-k&\leq d(p_m,p)\leq d(p_m,q_n)+d(q_n,q)+d(q,p)\\
&\leq \delta'+1+d_A(1,q)-n+2k+1,
\end{align*}
thus, $n-m\leq \delta'+1+d_A(1,q)-d_A(1,p)+2k+1\leq \delta'+2k+3.$

Hence, we have that 
$$d(p_n,q_n)\leq d(p_n,p_m)+d(p_m,q_n)\leq |m-n|+k+\delta'+1\leq \delta'+2k+3+k+\delta'+1=3k+2\delta'+4.$$
\qed

\section{Bounded Reduction Property}
 \label{secbrp}
 
In this section, we will present three (equivalent) geometric versions of the Bounded Reduction Property (also known as the Bounded Cancellation Lemma) for hyperbolic groups. The bounded reduction property has proved itself to be a most useful tool in studying the dynamics of (virtually) free group endomorphisms (see, for example, \cite{[BH92]}, \cite{[Coo87]},\cite{[GJLL98]},\cite{[Sil10]} \cite{[Sil13]}). We will later use it to prove the main result of this paper.

For free groups, the \emph{bounded cancellation lemma} is said to hold for an endomorphism $\varphi\in \End(F_n)$ if there is some constant $N\in \N$ such that $$|u^{-1}\wedge v|=0\Rightarrow |u^{-1}\varphi\wedge v\varphi|<N$$
holds for all $u,v\in F_n.$

Since, for $u,v\in F_n$, the Gromov product $(u|v)$ coincides with $u\wedge v$, a natural generalization for the hyperbolic case is
\begin{align*}
\exists N\in \N \, ((u\mid v)=0 \Rightarrow (u\varphi \mid v\varphi)\leq N).
\end{align*}

In Proposition 15, Chapter 5 in \cite{[GH90]}, it is shown that if $\varphi$ is a $(\lambda,K)$-quasi-isometric embedding, then there exists a constant $A$ depending on $\lambda, K$ and $\delta$, with the property that 
\begin{align}
\label{ghqie}
\frac 1 \lambda (u|v)-A\leq (u\varphi|v\varphi)\leq \lambda (u|v) +A.
\end{align}
So, if $H$ is a hyperbolic group and $\varphi:H\to H$ is a quasi-isometric embedding, then for every $p\geq 0$, there exists $q\geq 0$ so that 
\begin{align*}
(u\mid v)\leq p \Rightarrow (u\varphi \mid v\varphi)\leq q.
\end{align*}

We will now present a geometric formulation of the inequality $(u|v)\leq p$ .

\begin{lemma}
\label{gromovgeo}
Let $u,v\in H$ and $p\in \N$. Then the following are equivalent:
\begin{enumerate}[(i)]
\item $(u|v)\leq p$
\item for any geodesics $\alpha$ and $\beta$ from $1$ to $u^{-1}$ and $v$, respectively, we have that the concatenation
$$\xymatrix{
1 \ar[rr]^{\alpha} && u^{-1} \ar[rr]^{\beta} && u^{-1}v
}$$
is a $(1,2p)$-quasi-geodesic
\item there are geodesics $\alpha$ and $\beta$ from $1$ to $u^{-1}$ and $v$, respectively, such that the concatenation
$$\xymatrix{
1 \ar[rr]^{\alpha} && u^{-1} \ar[rr]^{\beta} && u^{-1}v
}$$
is a $(1,2p)$-quasi-geodesic
\end{enumerate}
\end{lemma}
\noindent\textit{Proof.}  It is clear that $(ii)\Rightarrow (iii)$. Now we prove that $(i)\Rightarrow (ii)$. Let $u,v\in H$ and $p\in \N$. Suppose that $(u|v)\leq p$ and take geodesics $\alpha$ and $\beta$ from $1$ to $u^{-1}$ and $v$, respectively. Take the concatenation $\zeta:[0,|u|+|v|]\to H$, where $0\zeta=1$, $(|u|)\zeta=u^{-1}$ and $(|u|+|v|)\zeta=u^{-1}v.$ 

 We will prove that $\zeta$ is a $(1,2p)$-quasi-geodesic, i.e., for $0\leq i\leq j\leq  |u|+|v|$, we have that $$j-i-2p\leq d(i\zeta,j\zeta)\leq j-i+2p.$$

We can assume that $0\leq i\leq |u|\leq j\leq |u|+|v|$. Clearly, we have that $ d(i\zeta,j\zeta)\leq j-i \leq j-i+2p$.

Suppose that $|j\zeta|<j-2(u|v)$ and consider a geodesic $\gamma:[0,|j\zeta|]$ from $1$ to $j\zeta$ and concatenate it with $\zeta\lvert_{[j,|u|+|v|]}$, which gives us a path of length 
\begin{align*}
|u|+|v|-j+|j\zeta|&<|u|+|v|-j+j-2(u|v)\\
&=|u|+|v|-2(u|v)\\
&=|u|+|v|-|u|-|v|+|u^{-1}v|\\
&=|u^{-1}v|
\end{align*}
 from $1$ to $u^{-1}v$ and that contradicts the definition of $d$. Thus, $|j\zeta|\geq j-2(u|v)$.

So, 
\begin{align*}
&d(u^{-1},j\zeta)+|u^{-1}v|-|v|-d(i\zeta,j\zeta)-d(1,i\zeta)
\\\leq\quad &d(u^{-1},j\zeta)+|u^{-1}v|-|v| -|j\zeta| 
\\=\quad &j-|u|+d(u,v)-|v| -|j\zeta|
\\=\quad &j -2(u|v)-|j\zeta|
\\\leq\quad & 0,
\end{align*}
thus, 
\begin{align*}
d(i\zeta,j\zeta)&\geq d(u^{-1},j\zeta)+|u^{-1}v|-|v|-d(1,i\zeta)\\
&=j-|u|+|u^{-1}v|-|v|-i\\
&=j-i-2(u|v)\\
&\geq j-i-2p
\end{align*}

To prove that $(iii)\Rightarrow (i)$, take geodesics $\alpha$ and $\beta$ from $1$ to $u^{-1}$ and $v$, respectively, such that the concatenation
$$\xymatrix{
1 \ar[rr]^{\alpha} && u^{-1} \ar[rr]^{\beta} && u^{-1}v 
}$$
is a $(1,2p)$-quasi-geodesic. Then $$d(u,v)=d(1,u^{-1}v)=d(0(\alpha+\beta),(|u|+|v|)(\alpha+\beta))\geq |u|+|v|-2p,$$ so $2(u|v)=|u|+|v|-d(u,v)\leq 2p$ and we are done.
\qed

 Let $\varphi:H\to H$ be a map. We say that \emph{the BRP holds for $\varphi$} if, for every $p\geq 0$ there is some $q\geq 0$ such that:
given two geodesics $u$ and $v$ such that $$\xymatrix{
1 \ar[rr]^u && u \ar[rr]^v && uv 
}$$ is a $(1,p)$-quasi-geodesic,  we have that given any two geodesics $\alpha$, $\beta$, from $1$ to $u\varphi$ and from $u\varphi$ to $(uv)\varphi$, respectively, the path 
$$\xymatrix{
1 \ar[rr]^{\alpha} && u\varphi \ar[rr]^{\beta} && (uv)\varphi 
}$$
is a $(1,q)$-quasi-geodesic.

\begin{proposition}
\label{geom1}
Let $H$ be a hyperbolic group and $\varphi:H\to H$ be a mapping. Then, the following are equivalent:
\begin{enumerate}[(i)]
\item BRP holds for $\varphi$
\item for every $p\geq 0$ there is some $q\geq 0$ such that, for all $u,v\in H$, we have that 
\begin{align}
\label{brpgromov}
(u\mid v)\leq p \Rightarrow (u\varphi \mid v\varphi)\leq q.
\end{align}
\end{enumerate}
\end{proposition}

\noindent\textit{Proof.} Let $H$ be a hyperbolic group, $\varphi:H\to H$ be a mapping and $p$ be a nonnegative integer. Suppose that the BRP holds for $\varphi$ and take $u,v\in H$ such that $(u|v)\leq p$. Then, by Lemma \ref{gromovgeo} and the BRP, there is some $q\geq 0$ such that, given geodesics $\alpha$ and $\beta$ from $1$ to $u^{-1}\varphi$ and from $u^{-1}\varphi$ to $(u^{-1}v)\varphi$, respectively, the concatenation $\zeta:[0,|u^{-1}\varphi|+|v\varphi|]\to H$ is a $(1,2q)$-quasi-geodesic. Using Lemma \ref{gromovgeo} again, we have that $(u\varphi|v\varphi)\leq q$.

Now, suppose that for every $p\geq 0$, there is some $q\geq 0$ such that (\ref{brpgromov}) holds and take geodesics $\alpha$, $\beta$ such that the concatenation
$$\xymatrix{
1 \ar[rr]^{\alpha} && u \ar[rr]^{\beta} && uv 
}$$
is a $(1,p)$-quasi-geodesic. In particular, it is also a $(1,2p)$-quasi-geodesic.  Then by Lemma \ref{gromovgeo}, we have that $(u^{-1}|v)\leq p$, so, using (\ref{brpgromov}), we have that $(u^{-1}\varphi|v\varphi)\leq q$, so by Lemma \ref{gromovgeo}, we have that the path
$$\xymatrix{
1 \ar[rr]^{\alpha} && u\varphi \ar[rr]^{\beta} && (uv)\varphi 
}$$
is a $(1,2q)$-quasi-geodesic for every geodesics $\alpha$, $\beta$ as above.
\qed

The following proposition is an immediate consequence of (\ref{brpgromov}) and Proposition \ref{geom1}.

\begin{proposition}
If $\varphi: H\to H$ is a quasi-isometric embedding, then the BRP holds for $\varphi.$ \qed
\end{proposition}

The next proposition shows that the bounded reduction property can be reduced to the case where $p=0$.

\begin{proposition}
\label{0chega}
Let $H$ be a hyperbolic group and $\varphi\in End(H)$. If the BRP holds for $\varphi$ for $p=0$, then it holds for every $p\in \N$.
\end{proposition}
\noindent\textit{Proof.} Let $H$ be a hyperbolic group and $\varphi\in End(H)$ and assume that the BRP holds when $p=0$. Let $p\in \N$, $u,v\in H$, and take two geodesics $\alpha$ and $\beta$ from $1$ to $u$ and from $u$ to $uv$, respectively, so that the concatenation $\alpha+\beta$ is a $(1,p)$-quasi-geodesic. By Proposition \ref{ftrav}, there is some $N\in \N$ such that for a $(1,p)$-quasi-geodesic $\xi$ starting in 1 and ending in $uv$, we have that 
\begin{align}
\label{ftp}
d(\xi^{[n]}, (\alpha+\beta)^{[n]})<N,
\end{align}
for every $n\in \N$.
We will prove that, there is some $q\in \N$ such that, given any two geodesics $\xi_1$ and $\xi_2$ from $1$ to $u\varphi$ and from $u\varphi$ to $(uv)\varphi,$ respectively, their concatenation is a $(1,q)$-quasi-geodesic. 

Take $\gamma$ to be a geodesic from $1$ to $uv$ and consider the factorization 
$$\xymatrix{
1 \ar[rr]^{\gamma^{[|u|]}} && x \ar[rr]^{\gamma_2} && uv .
}$$
Notice that both $\gamma^{[|u|]}$ and $\gamma_2$ are geodesics and their concatenation, $\gamma$ is also a geodesic.
In particular $\gamma$ is a $(1,p)$-quasi-geodesic with the same starting and ending points as the concatenation of $\alpha$ and $\beta$.  So, by, (\ref{ftp}) we have that $d(x,u)<N$.

Set $B_\varphi=\max\{|a\varphi|\,\big\lvert\,a\in \widetilde A\}$. Then, we have that $d(x\varphi,u\varphi)\leq B_\varphi d(x,u)<B_\varphi N.$ Let $\zeta_1$ and $\zeta_2$ be geodesics from $1$ to $x\varphi$ and from $x\varphi$ to $(uv)\varphi$, respectively. Since the BRP holds for $\varphi$ when $p=0$, we have that there is some constant $p_0$ such that the concatenation $\zeta_1+\zeta_2$ is a $(1,p_0)$-quasi-geodesic and that constant is independent from the choice of $\gamma$. 
So, we have that
\begin{align}
\label{qg1}
d(1,(uv)\varphi)\geq d(1,x\varphi)+d(x\varphi,(uv)\varphi)-p_0.
\end{align}
Since $d(u\varphi,(uv)\varphi)\leq d(u\varphi,x\varphi)+d(x\varphi,(uv)\varphi)\leq NB_\varphi+d(x\varphi,(uv)\varphi)$, we have that 
\begin{align}
\label{dt1}
d(x\varphi,(uv)\varphi)\geq d(u\varphi,(uv)\varphi)-NB_\varphi.
\end{align}
Similarly, we have that $d(1,u\varphi)\leq d(1,x\varphi)+d(x\varphi,u\varphi)\leq d(1,x\varphi)+NB_\varphi$, and so
\begin{align}
\label{dt2}
d(1,x\varphi)\geq d(1,u\varphi)-NB_\varphi.
\end{align}

Combining (\ref{qg1}) with (\ref{dt1}) and (\ref{dt2}), we have that 
\begin{align*}
d(1,(uv)\varphi)&\geq d(1,u\varphi)-NB_\varphi+d(u\varphi,(uv)\varphi)-NB_\varphi-p_0
\\&=d(1,u\varphi)+d(1,v\varphi)-2NB_\varphi-p_0.
\end{align*}
Hence 
\begin{align*}
2\left((u\varphi)^{-1}|v\varphi\right)&=d(1,u\varphi)+d(1,v\varphi)-d((u\varphi)^{-1},v\varphi)
\\&=d(1,u\varphi)+d(1,v\varphi)-d(1,(uv)\varphi)
\\&\leq 2NB_\varphi+p_0.
\end{align*}
By Lemma \ref{gromovgeo}, we have that the concatenation $\xi_1+\xi_2$ is a $(1,2NB_\varphi+p_0)$-quasi-geodesic.
\qed

If $(X,d)$ is $\delta$-hyperbolic and $\lambda \geq 1$, $K \geq 0$, it
follows from  \cite[Thm 1.7, Section III.H.3]{[BH99]}
that there exists a constant $R(\delta,\lambda,K)$,
depending only on $\delta,\lambda,K$, such that any geodesic
and $(\lambda,K)$-quasi-geodesic in $X$ having the same initial and
terminal points 
lie at Hausdorff distance $\leq R(\delta,\lambda,K)$ from each
other. This constant will be used in the proof of the next result.

We recall that for a geodesic $\alpha:[0,n]\to H$, we will often denote  its image by $\alpha$ as well.  We are now ready to present two more (equivalent) geometric formulations of the BRP.

\begin{theorem}
\label{geom2}
Let $\varphi\in \End(H)$. The following conditions are equivalent:
\begin{enumerate}[(i)]
\item the BRP holds for $\varphi$.
\item there is some $N\in \N$ such that, for all $x,y\in H$ and every geodesic $\alpha=[x,y]$, we have that $\alpha\varphi$ is at bounded Hausdorff distance to every geodesic $[x\varphi,y\varphi]$.
\item there is some $N\in \N$ such that, for all $x,y\in H$ and every geodesic $\alpha=[x,y]$, we have that $\alpha\varphi\subseteq \mc V_N(\xi)$ for every geodesic $\xi=[x\varphi,y\varphi]$.
 \end{enumerate}
\end{theorem}

\noindent\textit{Proof.} Clearly (ii) $\Rightarrow$ (iii).

(i) $\Rightarrow$ (ii) :  
Let $x, y\in H$ and $N\in \N$ given by the BRP when $p=0$. Consider geodesics $\alpha=[x,y]$ and $\xi=[x\varphi,y\varphi]$. Let $u\in\xi$ and $k=d(x\varphi,u)$. So, clearly, $d(x\varphi, y\varphi)\geq k$. Put $B_\varphi=\max\{|a\varphi|\,\big\lvert\,a\in \widetilde A\}$. We may assume that $k>B_\varphi$ since, otherwise $d(u,\alpha\varphi)\leq B_\varphi$. 
Since $d(x\varphi,x\varphi)=0<k$ and $d(y\varphi,x\varphi)\geq k$, there is some $n_k>0$ such that $d(x\varphi,(x(\alpha^{[n_k]}\pi))\varphi)<k$ and  $d(x\varphi,(x(\alpha^{[n_k+1]}\pi))\varphi)\geq k$ (notice that $n_k>0$ since  $d(x\varphi,(x(\alpha^{[1]}\pi))\varphi)\leq B_\varphi< k$ ). Consider the following factorization of $\alpha$
$$\xymatrix{
x \ar[rr]^{\alpha^{[n_k]}} && x_k \ar[rr]^{\alpha_k} && y 
}.$$ 

Using the BRP, we have that, given geodesics $\beta, \gamma$ from $x\varphi$ to $x_k\varphi$ and from $x_k\varphi$ to $y\varphi$, respectively, the concatenation 
$$\xymatrix{
x\varphi \ar[rr]^{\beta} && x_k\varphi \ar[rr]^{\gamma} && y\varphi 
}$$ 
is a $(1,N)$-quasi-geodesic. Set $x_{k+1}=x(\alpha^{[n_k+1]}\pi)$. We know that $d(x_k\varphi,x_{k+1}\varphi)\leq B_\varphi$ and so  $$k\leq d(x\varphi,x_{k+1}\varphi)\leq d(x\varphi,x_k\varphi)+B_\varphi.$$
Thus, we have that $d(x_k\varphi,x\varphi)\geq k-B_\varphi.$

Let $z=(\beta+\gamma)^{[k]}.$ Notice that, since $d(x_k\varphi,x\varphi)<k$, then $z\in \gamma.$
$$
\begin{tikzcd}[sep=large]
&u \ar[dr,bend left,"\alpha_k"]&\\
   x\varphi \ar[r,swap,bend left,"\beta"] \ar[ur,bend left,"\alpha^{[k]}"] &x_k\varphi  \ar[r,swap,bend left,"\gamma"] & y\varphi
\end{tikzcd}
$$
 Using the fellow traveller property for $(1,N)$-quasi-geodesics, there is some constant $M$ depending only on $N$ and $\delta$ such that 
$d(u,z)<M$. Since $\gamma$ is a geodesic, then $d(x_k\varphi,z)=k-d(x\varphi,x_k\varphi)\leq k-(k-B_\varphi)=B_\varphi$. Thus, $$d(u,x_k\varphi)\leq  d(u,z)+d(z,x_k\varphi)<M+B_\varphi$$ 
and $u\in \mc V_{M+B_\varphi}(\alpha\varphi)$. Since $u$ is an arbitrary element of $\xi$ such that $d(u,\alpha\varphi) >B_\varphi$, we have that $\xi\subseteq \mc V_{M+B_\varphi}(\alpha\varphi)$

Now, let $v\in \alpha$. Since the BRP holds, then, taking any geodesics, $\beta',\gamma'$ from $x\varphi$ to $v\varphi$ and from $v\varphi$ to $y\varphi$, the concatenation $$\xymatrix{
x\varphi \ar[rr]^{\beta'} && v\varphi \ar[rr]^{\gamma'} && y\varphi 
}$$ 
is a $(1,N)$-quasi-geodesic, so Haus$((\beta'+\gamma'),\xi)<R(\delta,1,N)$. In particular, $d(v\varphi,\xi)\leq R(\delta,1,N)$. Since $v$ is arbitrary, we have that $\alpha\varphi\subseteq \mc V_{R(\delta,1,N)}(\xi)$.

Hence Haus$(\xi,\alpha\varphi)\leq \max\{R(\delta,1,N),M+B_\varphi\}$.

(iii)$\Rightarrow$ (i) Take $N$ such that (iii) holds and $u,v\in H$. Let $\alpha=[1,u]$, $\beta=[u,uv]$ be such that $\alpha+\beta$ is a geodesic and consider $\gamma=[1,u\varphi]$ and $\zeta=[u\varphi,(uv)\varphi]$. We want to prove that there is some $M\in \N$ such that $\gamma+\zeta$ is a $(1,M)$-quasi-geodesic and that suffices by Proposition \ref{0chega}. From (iii), we have that $\alpha\varphi\subseteq \mc V_N(\gamma)$, $\beta\varphi\subseteq \mc V_N(\zeta)$ and $(\alpha+\beta)\varphi\subseteq \mc V_N(\gamma+\zeta)$. Since $\alpha+\beta$ is a geodesic, then, given a geodesic $\xi=[1,(uv)\varphi]$, we have that $((\alpha+\beta)\varphi)\subseteq \mc V_N(\xi)$ too. 

$$
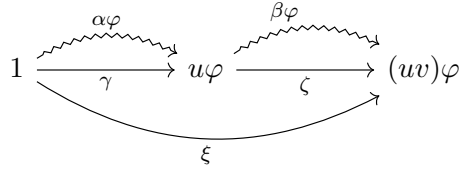
\begin{tikzcd}[sep=huge]
   1 \ar[r,bend left,squiggly,"\alpha\varphi"] \ar[r,swap,"\gamma"]\ar[rr,swap,bend right,"\xi"] &u\varphi  \ar[r,swap,"\zeta"] \ar[r,bend left,squiggly,"\beta\varphi"] & (uv)\varphi
\end{tikzcd}
$$

Now, we have that there is some $x_u\in\xi$ such that $d(u\varphi,x_u)< N$. Suppose that $d(1,x_u)\geq |\gamma|$ and denote $\xi^{[|\gamma|]}$ by $y$. So, $\gamma$ and $\xi_u=\xi^{[|x_u|]}$ are two geodesics with the same starting point that end at bounded distance. By the fellow traveller property, there is some $K\in \N$ such that $d\left(u\varphi,y)\right)<K$. So, we have that, $$|\zeta|=d(u\varphi,(uv)\varphi)\leq d(u\varphi,y)+d(y,(uv)\varphi)\leq K+d(y,(uv)\varphi).$$ Now, $$|\xi|=d(1,y)+d(y,(uv)\varphi)=|\gamma|+d(y,(uv)\varphi)\geq |\gamma|+|\zeta|-K.$$
If $d(1,x_u)<|\gamma|$, the same inequality can be obtained analogously, considering the geodesics $\zeta^{-1}$ and $\xi^{-1}$, since $|\gamma|+|\zeta|\geq |\xi|.$

So, we have that 
\begin{align*}
(u^{-1}\varphi|v\varphi)&=\frac 12 (d(1,u\varphi)+d(1,v\varphi)-d(u^{-1}\varphi,v\varphi))
\\&=\frac 12 (|\gamma|+|\zeta|-d(1,(uv)\varphi))=\frac 12 (|\gamma|+|\zeta|-|\xi|) \\
&\leq \frac K 2.
\end{align*}

By Lemma \ref{gromovgeo}, we have that $\gamma+\zeta$ is a $(1,K)$-quasi-geodesic.
\qed\\

A metric space $(X,d)$ is said to be a \emph{median space} if, for all $x,y,z\in X$, there is some 
unique point $\mu(x,y,z)\in X$, known as the \emph{median} of $x,y,z$, 
such that $d(x,y)=d(x,\mu(x,y,z))+d(\mu(x,y,z),y);$ $d(y,z)=d(y,\mu(x,y,z))+d(\mu(x,y,z),z);$ and $d(z,x)=d(z,\mu(x,y,z))+d(\mu(x,y,z),x)$. We call $\mu:X^3\to X$ the \emph{median operator} of the median space $X$. 

Coarse median spaces were introduced by Bowditch in \cite{[Bow13]}. Following the equivalent definition given in \cite{[NWZ19]}, we say that, given a metric space $X$, a \emph{coarse median on $X$} is a ternary operation $\mu:X^3\to X$ satisfying the following: 

there exists a constant $C\geq 0$ such that, for all $a,b,c,x\in X$, we have that 
\begin{enumerate}
\item $\mu(a,a,b)=a \text{ and } \mu(a,b,c)=\mu(b,c,a)=\mu(b,a,c);$
\item $d(\mu(\mu(a,x,b),x,c),\mu(a,x,\mu(b,x,c)))\leq C;$
\item $d(\mu(a,b,c),\mu(x,b,c))\leq Cd(a,x)+C.$
\end{enumerate} 

Given a group $G$, a \emph{word metric} on $G$ measures the distance of the shortest path in the Cayley graph of $G$ with respect to some set of generators, i.e., for two elements $g,h\in G$, we have that $d(g,h)$ is the length of the shortest word whose letters come from the generating set representing $g^{-1}h$.
Following the definitions in \cite{[Fio21]}, two coarse medians $\mu_1,\mu_2:X^3\to X$ are said to be at \emph{bounded distance} if there exists some constant $C$ such that $d(\mu_1(x,y,z),\mu_2(x,y,z))\leq C$ for all $x,y,z\in X$, and a \emph{coarse median structure} on $X$ is an equivalence class $[\mu]$ of coarse medians pairwise at bounded distance. When $X$ is a metric space and $[\mu]$ is a coarse median structure on $X$, we say that $(X,[\mu])$ is a \emph{coarse median space.} Following Fioravanti's definition, a \emph{coarse median group} is a pair $(G,[\mu])$, where $G$ is a finitely generated group with a word metric $d$ and $[\mu]$ is a $G$-invariant coarse median structure on $G$, meaning that for each $g\in G$, there is a constant $C(g)$ such that $d(g\mu(g_1,g_2,g_3),\mu(gg_1,gg_2,gg_3)\leq C(g)$, for all $g_1,g_2,g_3\in G$. The author in \cite{[Fio21]} also remarks that this definition is stronger than the original definition from \cite{[Bow13]}, that did not require $G$-invariance. Despite being better suited for this work, it is not quasi-isometry-invariant nor commensurability-invariant, unlike Bowditch's version.

An equivalent definition of hyperbolicity is given by the existence of a center of geodesic triangles (see, for example, \cite{[Bow91]}).
\begin{lemma}
\label{hyper2}
A group $G$ is hyperbolic if and only if there is some constant $K\geq 0$ for which every geodesic triangle has a $K$-center, i.e., a point that, up to a bounded distance, depends only on the vertices, and is $K$-close to every edge of the triangle. 
\end{lemma}
Given three points, the operator that associates the three points to the $K$-center of a geodesic triangle they define is coarse median. In fact, by Theorem 4.2. in \cite{[NWZ19]} it is the only coarse-median structure that we can endow $X$ with.

We will now show that the BRP coincides with coarse-median preservation. 

\begin{theorem}
\label{cmphyper}
Let $G$ be a hyperbolic group and $\varphi\in \End(G)$. Then the BRP holds for $\varphi$ if and only if $\varphi$ is coarse-median preserving.
\end{theorem}

\noindent \textit{Proof.}
Let $H=\langle A\rangle$ be a hyperbolic group and $\varphi\in \End(H)$. Suppose that the BRP holds for $\varphi$ and take $N$ given by version $(ii)$ of the BRP given by Theorem \ref{geom2}, $K$ be the constant given by Lemma \ref{hyper2}, $B_\varphi=\max\{d_A(1,a\varphi)\mid a\in A\}$,  $x,y,z\in G$ and geodesics $\alpha=[x,y]$, $\beta=[y,z]$, $\gamma=[z,x]$. Then by definition of the coarse-median operator, $d(\mu(x,y,z),\alpha)\leq K$, and so there is some $x'\in \alpha$ such that $d(x',\mu(x,y,z))\leq K$ and so 
\begin{align}
\label{median1}
d(x'\varphi, \mu(x,y,z)\varphi)\leq KB_\varphi.
\end{align} Now consider the triangle defined by $\alpha'=[x\varphi,y\varphi], \beta'=[y\varphi,z\varphi]$ and $\gamma'=[z\varphi,x\varphi]$. By the BRP, we have that $\alpha\varphi$ is at a Hausdorff distance at most $N$ from $\alpha'$, so there is some $x''\in \alpha'$ such that 
\begin{align}
\label{median2}
d(x'',x'\varphi)\leq N.
\end{align}

Combining (\ref{median1}) and (\ref{median2}), we get that $d(\mu(x,y,z)\varphi,x'')\leq KB_\varphi+N$, and so $d(\mu(x,y,z)\varphi,\alpha')\leq KB_\varphi+N.$ Similarly, we can prove that $d(\mu(x,y,z)\varphi,\beta')\leq KB_\varphi+N$ and that $d(\mu(x,y,z)\varphi,\gamma')\leq KB_\varphi+N$ and so $\mu(x,y,z)\varphi$ is a $(KB_\varphi+N)$-center of the triangle defined by $\alpha',\beta'$ and $\gamma'$. Since $\mu(x\varphi,y\varphi,z\varphi)$ is a $K$-center of the triangle, it is also a $(KB_\varphi+N)$-center and by Lemma 3.1.5 in \cite{[Bow91]}, we have that the distance between  $\mu(x\varphi,y\varphi,z\varphi)$  and  $\mu(x,y,z)\varphi$ is bounded and so $\varphi$ is coarse-median preserving.
 
Now, suppose that $\varphi$ is coarse-median preserving with constant $C\geq 0$. Take two points $x,y\in H$ and consider a geodesic $\alpha=[x,y].$ Let $x_i\in \alpha$. Then $\mu(x,x_i,y)=x_i$. Now consider a geodesic triangle given by $\alpha'=[x\varphi, x_i\varphi], \beta'=[x_i\varphi,y\varphi]$ and $\gamma'=[x\varphi,y\varphi]$. Since $\varphi$ is coarse-median preserving, we have that 
\begin{align*}
d(x_i\varphi,\gamma')\leq\; & d(x_i\varphi,\mu(x\varphi,x_i\varphi,y\varphi))+d(\mu(x\varphi,x_i\varphi,y\varphi),\gamma')\\
\leq\;&d(\mu(x,x_i,y)\varphi,\mu(x\varphi,x_i\varphi,y\varphi))+K\\
\leq\;& C+K.
\end{align*} 
Since $x_i\in \alpha$ is arbitrary, this means that $\alpha\varphi\subseteq \mc V_{C+K}(\gamma')$, which means that the BRP holds for $\varphi$, by Theorem \ref{geom2}.
\qed

So, combining Propositions \ref{geom1} and \ref{0chega} with Theorems \ref{geom2} and \ref{cmphyper}, we have proved the following result:
\begin{theorem}
Let $\varphi\in \End(H)$. The following conditions are equivalent:
\label{equivsbrp}
\begin{enumerate}[i.]
\item The BRP holds for $\varphi$.
\item The BRP holds for $\varphi$ when $p=0$.
\item $\forall\, p>0\,\exists\, q>0\,\forall\, u,v\in H\, ((u|v)\leq p\Rightarrow (u\varphi|v\varphi)\leq q)$.
\item $\exists\, q>0\,\forall\, u,v\in H\, ((u|v)=0\Rightarrow (u\varphi|v\varphi)\leq q)$.
\item there is some $N\in \N$ such that, for all $x,y\in H$ and every geodesic $\alpha=[x,y]$, we have that $\alpha\varphi$ is at bounded Hausdorff distance to every geodesic $[x\varphi,y\varphi]$.
\item there is some $N\in \N$ such that, for all $x,y\in H$ and every geodesic $\alpha=[x,y]$, we have that $\alpha\varphi\subseteq \mc V_N(\xi)$ for every geodesic $\xi=[x\varphi,y\varphi]$.
\item $\varphi$ is coarse-median preserving.
\end{enumerate}
\end{theorem}

\section{Uniformly continuous endomorphisms}
 \label{secuc}

We will start by describing when an endomorphism of a hyperbolic group admits a continuous extension to the completion.
It is well know by a general topology result \cite[Section XIV.6]{[Dug78]} that every uniformly continuous mapping $\varphi:H\to H'$ admits a continuous extension to the completion. So, we have the following lemma:

\begin{lemma}
Let $\varphi:H\to H'$ be a mapping of hyperbolic groups, and let $d$ and $d'$ be visual metrics on $G$ and $G'$, respectively. Then, the following conditions are equivalent:
\begin{enumerate}
\item $\varphi$ is uniformly continuous with respect to $d$ and $d'$;
\item $\varphi$ admits a continuous extension $\hat\varphi:(\widehat H,\hat d)\to (\widehat H',\hat d')$.
\end{enumerate}
\end{lemma}

A mapping $\varphi:(X,d)\to (X',d')$ between metric spaces satisfies a \emph{H\"older condition} of exponent $r>0$ if there exists a constant $K>0$ such that $$d'(x\varphi,y\varphi)\leq K(d(x,y))^r$$
for all $x,y\in X$. It clearly implies uniform continuity. We will show that in case of hyperbolic groups, the converse also holds.

In \cite{[AS16]}, the authors thoroughly study endomorphisms of hyperbolic groups satisfying a H\"older condition. In particular, they find several properties equivalent to satisfying a H\"older condition. An endomorphism $\varphi$ of $H$ is \emph{virtually injective} if its kernel is finite.

\begin{theorem}\cite[Thm 4.3]{[AS16]}
\label{asequivs}
Let $\varphi$ be a nontrivial endomorphism of a hyperbolic group $G$ and
let $d \in V^A(p,\gamma,T)$ be a visual metric on $G$. Then the 
following conditions are equivalent:
\begin{enumerate}
\item[(i)] $\varphi$ satisfies a H\"older condition with respect to $d$;
\item[(ii)] $\varphi$ admits an extension to $\widehat{G}$ satisfying
  a H\"older condition with respect to $\widehat{d}$; 
\item[(iii)] there exist constants $P > 0$ and $Q \in \R$ such that
$$
P(g\varphi|h\varphi)_p^A + Q \geq (g|h)_p^A
$$
for all $g,h \in G$;
\item[(iv)] $\varphi$ is a quasi-isometric embedding of $(G,d_A)$ into itself;
\item[(v)] $\varphi$ is virtually injective and $G\varphi$ is a quasiconvex
  subgroup of $G$.
\end{enumerate}
\end{theorem}

The authors in \cite{[AS16]} conjecture that every uniformly continuous endomorphism satisfies a H\"older condition. We will give a positive answer to that problem later in this section.

We now present a natural result following from Theorem \ref{geom2}.
\begin{corollary}
\label{brpqc}
Let $\varphi\in \End(H)$ such that the BRP holds for $\varphi$. Then $H\varphi$ is quasiconvex.
\end{corollary}
\noindent\textit{Proof.} Let $x, y\in H$ and take $N\in \N$ given by condition (ii) of Theorem \ref{geom2}. Consider geodesics $\alpha=[x,y]$ and $\xi=[x\varphi,y\varphi]$. We have that  Haus$(\xi,\alpha\varphi)\leq N$. Let $u\in\xi$. Then $$d(u,H\varphi)\leq d(u,\alpha\varphi)\leq \text{Haus}(\xi,\alpha\varphi)\leq N.$$
\qed

Lemma 4.1 in \cite{[AS16]} states that a uniformly continuous endomorphism is virtually injective. So, next we will prove that uniform continuity implies the BRP. In that case, it follows that  uniformly continuous endomorphisms are precisely the ones satisfying a  H\"older condition.

Let $p,q,x,y\in H$. We have that 
\begin{align*}
(x|y)_p&=\frac 12 (d_A(p,x)+d_A(p,y)-d_A(x,y))\\
&\leq \frac 12  (d_A(q,x)+d_A(q,y)-d_A(x,y)+2d_A(p,q))\\
&=(x|y)_q+d_A(p,q)
\end{align*}
Similarly, $(x|y)_q\leq (x|y)_p +d_A(p,q)$, so 
\begin{align}
\label{trocaponto}
(x|y)_q-d_A(p,q)\leq (x|y)_p\leq (x|y)_q+d_A(p,q)
\end{align}

\begin{proposition}
\label{uninj}
Let $G=\langle A\rangle$ and $H=\langle B\rangle$ be hyperbolic groups and consider visual metrics $d_1\in V^A(p,\gamma,T)$ and $d_2\in V^B(p',\gamma',T')$ on $G$ and $H$, respectively. Let $\varphi: (G,d_1)\to (H,d_2)$ be an injective uniformly continuous homomorphism. Then, for every $M\geq 0$, there is some $N\geq 0$ such that $$(u|v)^A\leq M \Rightarrow (u\varphi | v\varphi)^B\leq N$$
holds for every $u,v\in G$.
\end{proposition}
\noindent\textit{Proof.}
Since $\varphi$ is uniformly continuous, by a general topology result it admits a continuous extension $\hat\varphi:(\hat G,\hat d_1)\to(\hat H,\hat d_2)$. Since $\hat \varphi$ is a continuous map between compact spaces, then it is closed, and so it has a closed (thus compact) image.

Now, restricting the codomain of $\hat\varphi$ to the image, we have a continuous bijection between compact spaces, and so it is a homeomorphism. Its inverse, $\psi: \text{Im}(\hat\varphi)\to \hat G$ is a continuous map between compact spaces, hence uniformly continuous. So, the restriction $\psi':(\text{Im(}\varphi\text{)},d_2)\to  (G,d_1)$ is also uniformly continuous, i.e., 
\begin{align*}
\forall \varepsilon >0\,\exists \delta>0 \, (d_2(x,y)<\delta\Rightarrow d_1(x\psi',y\psi')<\varepsilon),
\end{align*}
which, by construction of $\psi'$, means that 
\begin{align}
\label{uc1}
\forall \varepsilon >0\,\exists \delta>0 \, (d_2(x\varphi,y\varphi)<\delta\Rightarrow d_1(x,y)<\varepsilon).
\end{align}
Using (\ref{ineq}), we have that (\ref{uc1}) is equivalent to 
\begin{align}
\label{uc2}
\forall M\in \N\,\exists N\in \N \, ((x\varphi|y\varphi)_{p'}^B>N\Rightarrow (x|y)_p^A>M).
\end{align}
Since $p$ and $p'$ are fixed, we can change the basepoint to $1$ using (\ref{trocaponto}). So, (\ref{uc2}) becomes equivalent to 
\begin{align*}
\forall M\in \N\,\exists N\in \N \, ((x|y)^A\leq M\Rightarrow (x\varphi|y\varphi)^B\leq N ).
\end{align*}
\qed

\begin{proposition}
\label{ucbrp}
Let $d\in V^A(p,\gamma,T)$ be a visual metric on $H$ and let $\varphi$ be a uniformly continuous endomorphism of $H$ with respect to $d$. Then, the BRP holds for $\varphi$.
\end{proposition}
\noindent\textit{Proof.} If $\varphi$ is injective it follows from Proposition \ref{uninj}.
Now, in case $\varphi$ is not injective, by Lemma 4.1 in \cite{[AS16]}, it must have finite kernel $K$. Consider $\pi: H\to \faktor  HK$ to be the projection and the geodesic metric $d_{A\pi}$ on the quotient. Let $\varphi': \left(\faktor{H}{K},d_{A\pi}\right)\to (H,d_A)$ be the injective homomorphism induced by $\varphi$.

$$
\begin{tikzcd}[sep=large]
  H       \ar[r,"\varphi"] \ar[d,swap,"\pi"] & H                                 \\
   \faktor{H}{K}  \ar[ur,"\varphi'"] &    
\end{tikzcd}
$$

Let 
$$L = \max\{ d_A(1,x) \mid x \in K\}.$$
Let $g,h \in H$. We claim that
\begin{align}
\label{fq2}
d_A(g,h)-L \leq d_{A\pi}(g\pi,h\pi) \leq d_A(g,h).
\end{align}

Since $h = ga_1\ldots a_n$ implies $h\pi = (ga_1\ldots a_n)\pi$ for
all $a_1,\ldots, a_n \in \widetilde{A}$, we have $d_{A\pi}(g\pi,h\pi) \leq
d_A(g,h)$.

Write $h\pi = (gw)\pi$, where $w$ is a word on $\widetilde{A}$ of minimum
length. Then $h = gwx$ for some $x \in K$ and so
$$d_A(g,h) \leq d_A(g,gw) + d_A(gw,gwx) = 
d_A(1,w) + d_A(1,x) \leq |w| + L.$$
By minimality of $w$, we have actually $|w| = d_{A\pi}(g\pi,h\pi)$ and
thus (\ref{fq2}) holds. 

This means that $(H,d_A)$ and $\left(\faktor{H}{K},d_{A\pi}\right)$ are quasi-isometric. In particular, it yields that $\faktor{H}{K}$ is hyperbolic.

Now, take $d'\in V^{A\pi}(p'\pi,\gamma',T')$ to be a visual metric on $\faktor{H}{K}$. For every $u,v,p\in H,$ we have that 
\begin{align*}
(u\pi|v\pi)_{p\pi}^{A\pi}&=\frac 12 (d_{A\pi}(p\pi,u\pi)+d_{A\pi}(p\pi,v\pi)-d_{A\pi}(u\pi,v\pi))\\
&\leq \frac 12 (d_A(p,u)+d_A(p,v)-d_A(u,v)+L)\\
&=(u|v)_p^A+\frac L2
\end{align*}

From uniform continuity of $\varphi$ with respect to $d$, we get that 
\begin{align*}
\forall M\in \N\,\exists N\in \N \, ((x|y)_{p}^A>N\Rightarrow (x\varphi|y\varphi)_p^A>M).
\end{align*}
It follows that 
\begin{align*}
\forall M\in \N\,\exists N\in \N \, ((x\pi|y\pi)_{p\pi}^{A\pi}>N\Rightarrow (x\varphi|y\varphi)_p^A>M),
\end{align*}
and so $\varphi'$ is uniformly continuous with respect to $d'$ and $d$. It follows from Proposition \ref{uninj} that for every $p\geq 0$, there is some $q\geq 0$ such that 

\begin{align}
\label{qqq}
(u\pi|v\pi)^{A\pi}\leq p \Rightarrow (u\pi\varphi' | v\pi\varphi')^A\leq q
\end{align}
holds for all $u\pi,v\pi\in \faktor{H}{K}$.

Take $u,v\in H$ such that $(u|v)^A=0$. Then 
\begin{align*}
(u\pi|v\pi)^{A\pi}&\leq(u|v)^A+\frac L 2=\frac L2.
\end{align*}

So, by (\ref{qqq}), there is some $q$ which does not depend on $u,v$ such that $(u\varphi|v\varphi)^A\leq q$.
 By Proposition \ref{0chega}, the BRP holds for $\varphi$.

\qed

We can now answer Problem 6.1 left by the authors in \cite{[AS16]}.
\begin{theorem}
Let $d\in V^A(p,\gamma,T)$ be a visual metric on $H$ and let $\varphi$ be an endomorphism of $H$. Then $\varphi$ is uniformly continuous with respect to $d$ if and only if the conditions from Theorem \ref{asequivs} hold.
\end{theorem}
\noindent\textit{Proof.} It is straightforward to see that condition (i) from Theorem \ref{asequivs} implies uniform continuity. Now, if $\varphi$ is uniformly continuous, by  Lemma 4.1 in \cite{[AS16]}  it must be virtually injective and combining Proposition \ref{ucbrp} with Corollary \ref{brpqc}, we have that $H\varphi$ is quasiconvex, so condition (v) of Theorem \ref{asequivs} holds.
\qed

We now present a visual representation of these properties, where the shaded region represents the nontrivial uniformly continuous endomorphisms. Indeed, we have proved that every nontrivial uniformly continuous endomorphism satisfies the BRP (Proposition \ref{ucbrp}) and that every endomorphism satisfying the BRP must have quasiconvex image (Corollary \ref{brpqc}). Theorem \ref{cmphyper} establishes an equivalence between the BRP and coarse-median preservation (CMP). 
Also, every virtually injective endomorphism with quasiconvex image must be uniformly continuous by Theorem \ref{asequivs}. In \cite{[AS16]}, the authors give an example of an injective endomorphism of a torsion-free hyperbolic group with non quasiconvex image. So, unlike the case of virtually free groups, the BRP does not hold in general for injective endomorphisms of hyperbolic groups, not even when restricted to torsion-free hyperbolic groups. In the virtually free groups case, that does not happen, as every virtually injective endomorphism is uniformly continuous \cite{[Sil13]}.

Taking $\varphi:F_3\to F_3$ defined by $a\mapsto a$, $b\mapsto b$ and $c\mapsto 1$, we have that $F_3\varphi=\langle a,b\rangle$. Since $F_3\varphi$ is finitely generated and the standard embedding $\langle a,b\rangle \xhookrightarrow{} F_3$ is a quasi-isometric embedding, then $\varphi$ has quasiconvex image. But the BRP does not hold for $\varphi$ since $|cb^n\wedge b^n|=0$ and $|(cb^n)\varphi\wedge b^n\varphi|=|b^n\wedge b^n|=n$, which can be arbitrarily large. 

It is easy to find examples of endomorphisms for which the BRP holds that are not virtually injective, even for virtually free groups, by taking an endomorphism with finite image. For example, take $H=\Z\times \Z_2$ and $\varphi$ defined by $(n,0)\mapsto (0,0)$ and $(n,1)\mapsto (0,1)$. Then, the BRP holds for $\varphi$ and its kernel is infinite (in particular, it can't be uniformly continuous). 
\begin{figure}[H]
\begin{center}
\begin{tikzpicture}
	\begin{scope} [fill opacity = 0.4]

    \draw[ draw = black] (-4.5,-2) rectangle (1.5,4);
        \draw[ draw = black] (-2.2,-1) rectangle (1.5,3);

    \draw[ draw = black]  (1.5,1) ellipse (2 and 1);

   \node[opacity=1] at (-3.5,1) {\small\text{\parbox{1.7cm} {\centering Quasiconvex image}}};
    \node[opacity=1] at (-1.3,1) {\small\text{\parbox{1.7cm} {\centering BRP, CMP}}};

        \node[opacity=1] at (1.5,1) {\small\text{Virtually Injective}};

        \clip (-2.2,-1) rectangle (1.5,3);
    \fill[ pattern=north west lines]  (1.5,1) ellipse (2 and 1);
    \end{scope}

\end{tikzpicture}
\caption{Nontrivial uniformly continuous endomorphisms of hyperbolic groups}
\end{center}
\end{figure}
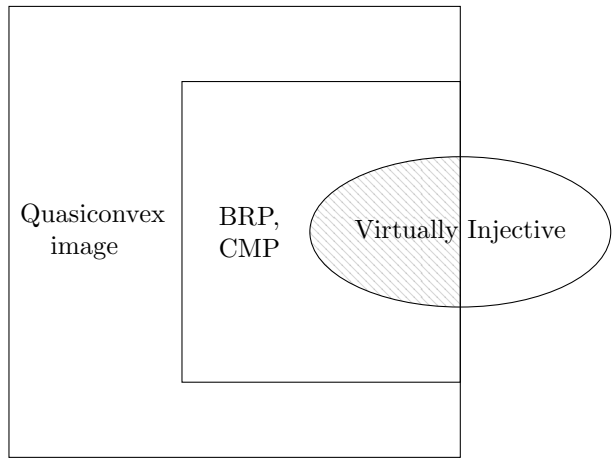

So, for hyperbolic groups, we have the figure above in which every region is nonempty. 
Notice that, for virtually free groups, the only difference is that nontrivial uniformly continuous endomorphisms are precisely the virtually injective ones and the BRP holds for all of them. 

In the case of free groups, it is even simpler as, for every nontrivial endomorphism, the properties of being injective, uniformly continuous and satisfying the BRP are equivalent. Indeed, it is well-known that injective endomorphisms coincide with uniformly continuous ones and that the BRP holds for this class. It is easy to see that the converse also holds. For $n\geq 2$, let $X=\{x_1,\ldots, x_n\}$ be a finite alphabet and $F_n=\langle X \rangle $ be a free group of rank $n$. If a nontrivial endomorphism $\varphi\in \End(F_n)$ is not injective, then there is some $w\in \text{Ker}(\varphi)$ such that $w\neq 1$ and some letter $a$  such that $a\varphi\neq1$. Let $p=|w|$. We have that $w$ is not a proper power of $a$ since, in that case we would have that $w\varphi$ would be a proper power of $a\varphi$ and so, nontrivial. So, we have that, for arbitrarily large $m\in \N$, $|wa^m\wedge a^m|<|w|$, but $|(wa^m)\varphi\wedge a^m \varphi|=|a^m\varphi|$, which is arbitrarily large. So, nontrivial endomorphisms for which the BRP holds  in a free group of finite rank are precisely the injective ones.

In \cite{[Pau89]}, Paulin proved that $\Fix(\varphi)$ is finitely generated if $\varphi\in \Aut(H)$. We remark that its proof also yields the result for quasi-isometric embeddings. So, we have proved the following result.
\begin{theorem} 
Let $\varphi\in \End(H)$ be an endomorphism admitting a continuous extension $\hat\varphi:\widehat H\to\widehat H$ to the completion of $H$. Then, $\Fix(\varphi)$ is finitely generated.\qed
\end{theorem}

\section*{Acknowledgements}

The author is grateful to  Armando Martino and to Pedro Silva for fruitful discussions of these topics, which greatly improved the paper.

The author was  supported by the grant SFRH/BD/145313/2019 funded by Funda\c c\~ao para a Ci\^encia e a Tecnologia (FCT).%

\end{document}